\documentclass[12pt]{article}
\usepackage{a4wide}

\usepackage{amssymb}
\usepackage{amsfonts}
\usepackage{amsmath}
\usepackage{amsthm}
\usepackage{graphicx}

\usepackage{empheq}

\begin{document}

\title{{\bfseries Time--Varying Capacitors as Lossless\\ Two--Port Devices}}

\author{Dimitri Jeltsema}

\date{Balanced Energy Systems\\ HAN University of Applied Sciences\\ P.O.~Box 2217, 6802 CE Arnhem, The Netherlands\\Email: d.jeltsema@han.nl}

\maketitle

\begin{abstract}
Inspired by recent assertions that imply incorrectness of the classical constitutive relationship of time--varying capacitance, energetic and system--theoretic perspectives are invoked that, instead of incorrectness, show incompleteness of the associated one--port device model. Based on dissipativity and passivity theory, an energetic analysis is presented that naturally leads to a complete two--port device model that is consistent with the laws of thermodynamics.
\end{abstract}

\section{Motivation}\label{sec:intro}

Recently, the paper ``Origin of the Curie--von Schweidler law and the fractional capacitor from time--varying capacitance''  
\cite{Vikash2022} asserted arguments to obtain a physical interpretation to this century--old law of dielectrics. The exposition in \cite{Vikash2022} starts by arguing  that the derivation of this law requires a revision of the classical charge--voltage relation of linear and time--varying capacitors, 
\begin{equation}\label{eq:QvsV}
Q(t) = C(t)V(t),
\end{equation}
where $Q$ represents the charge, $C$ the capacitance, $V$ the voltage, and $t$ represents time. The proposed revision of (\ref{eq:QvsV}) is expressed in \cite{Vikash2022} as a convolution between the time--varying capacitance and the time--derivative of the associated voltage, i.e., 
\begin{equation}\label{eq:Q*V}
Q(t) = C(t) * \dot{V}(t).
\end{equation}
A similar type of revision is proposed in \cite{Fouda2020}. 

To substantiate the arguments against (\ref{eq:QvsV}), the paper \cite{Vikash2022}, and additionally the response \cite{Vikash2022_response} to another critical analysis \cite{Allagui2022}, discusses a simple example in which an external voltage is applied to a linearly increasing capacitance of the form $C(t) = C_0 + \phi t$, with $C_0,\phi >0$. The resulting current is then obtained by taking the time--derivative of the charge, i.e., 
\begin{equation}\label{eq:LTVcap}
I(t) = \dot{Q}(t) = \frac{d}{dt}\big[C(t)V(t)\big]  = C(t)\dot{V}(t) + \dot{C}(t)V(t),
\end{equation}
which, for the assumed linearly increasing capacitance, implies that
\begin{equation}\label{eq:LTVcap_phi}
I(t) = C_0 \dot{V}(t) + \phi t \dot{V}(t) + \phi V(t).
\end{equation}

Besides the fact that (\ref{eq:LTVcap}) is advocated in classical textbooks like \cite{ChuaBook1987,Desoer,Valkenburg1955}, and is also used in simulation programs like Matlab, Micro--Cap, or Pspice \cite{Biolek}, the argumentation in \cite{Vikash2022} proceeds with the evaluation of (\ref{eq:LTVcap_phi}) at $t=0$ and points out an apparent inconsistency as the initial capacitance $C(0) = C_0$ suggests that the initial current should read as
\begin{equation}\label{eq:LTVcapI0_wrong}
I(0) = C_0\dot{V}(0),
\end{equation} while (\ref{eq:LTVcap}) dictates that
\begin{equation}\label{eq:LTVcap0}
I(0) = C_0\dot{V}(0) + \phi V(0).
\end{equation}
Based on these considerations, \cite{Vikash2022} labels the term $\phi V(t)$ in (\ref{eq:LTVcap_phi}) as an anomaly arising from the constitutive relationship (\ref{eq:QvsV}), and concludes (\ref{eq:QvsV}) to be erroneous.  

Although at first glance this observation might seem conspicuous and convincing, there is an elementary mathematical flaw in this reasoning that can be easily resolved by a proper classification of (\ref{eq:LTVcap}). Nevertheless, there are indeed some fundamental energetic and system--theoretic issues with the device model (\ref{eq:LTVcap}) that will be the main focus of the remainder of the present work. 

\section{Outline and contribution}

The contribution of this paper is multi--folded. In Section \ref{sec:math_flaw}, it will be argued that (\ref{eq:LTVcap}) naturally defines an elementary linear ordinary differential equation that generally consists of both a homogeneous and a particular solution---something that appears to be overlooked in \cite{Vikash2022}. This renders the initial arguments in \cite{Vikash2022} to advocate for (\ref{eq:Q*V}) invalid, and the discussion is proceeded with a short subsection about physical causality and concludes by showing that (\ref{eq:LTVcap}) violates gauge symmetry. Section \ref{sec:passivity} addresses the fact that (\ref{eq:LTVcap}) constitutes an electrical one--port device that is physically incomplete since the mechanism that causes a capacitor to change its capacitance is missing in the model. Based on the energetic and system--theoretic notions of dissipativity and passivity, it is shown that (\ref{eq:LTVcap}) violates the First Law of thermodynamics when the capacitance varies in a cyclic fashion. In Section \ref{sec:completeness}, the device model (\ref{eq:LTVcap}) is completed by considering a time--varying capacitor as a two--port device, with an electrical port and a port that provides the `mechanical' work to accomplish the change in the capacitance. Additionally, Section \ref{sec:2Cap_paradox} provides a consistent resolution to a long standing debate of a non--cyclic example that is (close to) equivalent with the peculiar so--called ``Two--Capacitor Paradox.''   
 The paper is concluded with some final remarks and a short discussion about the proposed convolution relation (\ref{eq:Q*V}). 

\section{First get the mathematics right\ldots}\label{sec:math_flaw}

\ldots the rest is physics.\footnote{Conversely, after Rudolf K\'alm\'an's famous quote: ``First get the physics right, the rest is mathematics.''} Indeed, it is not difficult to observe that (\ref{eq:LTVcap}) rather---and most naturally---defines a linear ordinary differential equation (ODE)  
\begin{equation}\label{eq:diff_VC}
\boxed{\dot{V}(t) + \frac{\dot{C}(t)}{C(t)}V(t) = f(t)}
\end{equation}
with forcing term $f(t) = I(t)/C(t)$. The solution to such standard ODE's can be found in virtually any undergraduate calculus textbook. Indeed, using the method of integrating factors, the general solution to (\ref{eq:diff_VC}) is given by (see e.g., \cite[Section 9.5]{Stewart})
\begin{equation*}
V(t) = e^{-\alpha(t)} \int e^{\alpha(t)}f(t)dt + \gamma e^{-\alpha(t)}, \ \dot{\alpha}(t) = \frac{\dot{C}(t)}{C(t)},
\end{equation*} 
where the constant $\gamma$ depends on the initial conditions. Clearly, this solution consists of two parts: a particular (forced) solution and a homogeneous (unforced) solution. The current $I(t)$ is the independent part of the forcing term $f(t)$ that determines, together with $V(0)$, the time--evolution of $V(t)$. Hence, in the voltage--current equation (\ref{eq:LTVcap}), the independent variable is the supplied current, whereas the dependent variable is the resulting capacitor voltage. To illustrate what goes wrong with the reasoning in \cite{Vikash2022}, let us return to the example of Section \ref{sec:intro}.

\subsection{The complete solution of (\ref{eq:LTVcap})}\label{sec:Viskash_example}

Continuing the example of \cite{Vikash2022} and \cite{Vikash2022_response}, let $V(t) = at$, with $a>0$, and $C(t) = C_0 + \phi t$ as before. Substitution of the latter into (\ref{eq:LTVcap}) suggest that the current will be given by
\begin{equation}\label{eq:Iwrong}
I(t) = (C_0 + 2\phi t)a.
\end{equation}
Along the lines of \cite{Vikash2022}, one would be tempted to conclude that $a = \dot{V}(t)$ and that the two terms on the right--hand side of (\ref{eq:Iwrong}) can be interpreted as two capacitive currents that flow through two capacitors that are arranged in a parallel combination of a constant capacitor $C_0$ and a time--varying capacitor $2\phi t$. Such interpretation would be conflicting with the fact that $C(t)$ equals $C_0 + \phi t$ instead of $C_0 + 2\phi t$ (as is suggested in Fig.~1 of \cite{Vikash2022}).  

However, this reasoning fails since the mapping $at \mapsto I(t)$ that results in (\ref{eq:Iwrong}) only considers the particular solution of the problem. Indeed, if we conversely substitute (\ref{eq:Iwrong}) into (\ref{eq:diff_VC}), then the \emph{complete} solution will be
\begin{equation}
V(t) = at + \frac{C_0}{C(t)}V(0) \ \Rightarrow \ \dot{V}(t) = a - \frac{C_0}{C^2(t)}\phi V(0),
\end{equation}
instead of just $\dot{V}(t) = a$. Hence, an equivalent parallel arrangement of two capacitors connected to a voltage source $V(t) = at$ breaks down, and so does the apparent contradiction between equations (\ref{eq:LTVcapI0_wrong}) and (\ref{eq:LTVcap0})  at $t=0$. The reason why the author of \cite{Vikash2022} was escorted to this apparent inconsistency in the first place can perhaps be explained in terms of (physical) causality. 

\subsection{Physical causality}

Although the time--differentiation in (\ref{eq:LTVcap}) of the voltage leading to the current is (in principle) a causal operation, it is not physically realizable as such. This becomes clear if one considers a (constant) capacitor directly in series with an ideal battery. It is not physically possible to charge a capacitor instantly. Indeed, whether the capacitance is constant or not, if a constant voltage is applied to (\ref{eq:LTVcap}), then the associated current is instantly infinite (i.e., the current would behave, at least theoretically, as a Dirac impulse). Moreover, by solely considering the constitutive relation between charge and voltage (\ref{eq:QvsV}), the dynamical aspect of charging is neglected. Indeed, a capacitor, either linear or nonlinear, time--varying or time--invariant, accumulates charge according to
\begin{equation}\label{eq:QvsI}
Q(t) = Q(0) + \int\limits_0^t I(\tau)d\tau, 
\end{equation}
which is not an instant process. 

\subsection{Violation of gauge symmetry}

Gauge symmetry \cite{Perkins1982} requires that the laws of physics are invariant under transformations like $V(t) \mapsto V(t) + \Psi$, where $\Psi$ is an arbitrary (constant) reference potential. This implies that no experiment should be able to measure the absolute potential, without reference to some external standard such as an electrical ground or datum point. Moreover, if gauge symmetry holds and energy is conserved, then charge must be conserved \cite{Perkins1982}. 

To test if (\ref{eq:LTVcap}) satisfies this principle, let us replace $V$ by $V+\Psi$ and evaluate the result, i.e.,
\begin{equation*}
I(t) = \frac{d}{dt}\big[C(t)(V(t) + \Psi)\big] =  C(t)\dot{V}(t) + \dot{C}(t)V(t) + \boxed{\dot{C}(t)\Psi}
\end{equation*}
which readily implies that
\begin{equation}\label{eq:gauge_fail}
\frac{d}{dt}\big[C(t)V(t)\big] \neq \frac{d}{dt}\big[C(t)\big(V(t) + \Psi\big)\big], 
\end{equation}
unless $C(t)$ is a constant! 

Apart from conservation of charge, the next question that arises is what this means for the energy conservation properties associated to the device model (\ref{eq:LTVcap}). 

\section{Violation of the First Law of thermodynamics}\label{sec:passivity}

For a device model to be physically consistent with the First Law of thermodynamics, we need to assure that the energy supplied to the device always equals, or is larger than, the amount of stored energy, i.e., the device must be \emph{passive}.\footnote{See the Appendix for a brief outline of dissipativity and passivity theory and some selected references about these concepts.} For the device model (\ref{eq:LTVcap}), or its equivalent (\ref{eq:diff_VC}), that is initially at rest, passivity amounts to the requirement that
\begin{equation}\label{eq:passivity}
\int\limits_{0}^{t} V(\tau)I(\tau) d\tau \geq 0,
\end{equation}
for all $t \geq 0$.  

\subsection{One--port formulation}

An equivalent way to formulate the current--voltage relationship (\ref{eq:LTVcap}) is by selecting the charge $Q(t)$ as an internal state variable, and define the linear and time--varying \emph{one--port} input--state--output system\footnote{This model is in the form of a so--called input--state--output \emph{port-Hamiltonian system}. The interested reader is referred to \cite{passivitybook} or \cite{NOW_book} for details of such system descriptions.}
\begin{equation}\label{eq:CpH}
\begin{aligned}
\dot{Q}(t) &= I(t) \ \ \text{(input)},\\
\text{(output)} \ \ V(t) &= \frac{\partial S}{\partial Q}\big(Q(t),t\big),
\end{aligned}
\end{equation}
with conjugated port--variables $(I,V)$ and 
\begin{equation}\label{eq:TVC_S}
S\big(Q(t),t\big) = \frac{1}{2C(t)}Q^2(t)
\end{equation}
represents the total energy stored in the capacitor at any instant of time $t$.
Indeed, differentiating (\ref{eq:TVC_S}) along the trajectories of (\ref{eq:CpH}) yields
\begin{equation}\label{eq:CpH_dotS_full}
\frac{d}{dt}S\big(Q(t),t\big) = V(t)I(t) - \frac{1}{2}\dot{C}(t)V^2(t). 
\end{equation}
The requirement for (\ref{eq:CpH}) to be passive is that 
\begin{equation}\label{eq:pass_ineq}
\frac{d}{dt}S\big(Q(t),t\big) \leq V(t)I(t),
\end{equation}
and lossless if the latter holds with equality. 

Obviously, for a constant capacitor $C(t) = C_0$, with $C_0>0$, this is trivially satisfied as (\ref{eq:pass_ineq}) holds with equality and therefore implies that constant capacitors are lossless devices (as should be expected). On the other hand, time--varying capacitors only satisfy (\ref{eq:pass_ineq}) as long as $C(t)$ is such that 
\begin{equation}\label{eq:dotC_cond}
\frac{1}{2}\dot{C}(t)V^2(t) \geq 0,
\end{equation}
for all time $t$. 

For the example of Section \ref{sec:Viskash_example}, this requirement is readily satisfied as $\phi V^2(t) \geq 0$, for all admissible capacitor voltages $V(t)$ and $\phi >0$. Hence, this hypothetical device is not purely lossless, but generally dissipates energy in moving from one state to another. The situation for the practically relevant class of linear capacitors that vary in a cyclic fashion, however, turns out to be quite different. 

\subsection{Linear cyclic capacitance}

Consider again the current--voltage equation (\ref{eq:LTVcap}) or, equivalently, the linear ODE (\ref{eq:diff_VC}) or the one-port formulation (\ref{eq:CpH}). Let us first consider a time--interval $[t_1,t_2]$ in which $C(t_1) = C(t_2)$ and $V(t_1) = V(t_2)$. Then, instead of starting from the total stored energy (\ref{eq:TVC_S}), consider the energy that is supplied to the device by the electrical port $(V,I)$.  Over one cycle the supplied energy is given by
\begin{equation}\label{eq:LTVcap_suppliedE}
\begin{aligned}
\oint\limits_{t_1}^{t_2} V(t)I(t) dt &= \oint\limits_{t_1}^{t_2} V(t)\frac{d}{dt}\big[C(t)V(t)\big]dt\\
&= C(t)V^2(t)\bigg|_{t_1}^{t_2} - \oint\limits_{t_1}^{t_2} C(t)V(t)\dot{V}(t)dt\\
&= \oint\limits_{t_1}^{t_2}\frac{1}{2} \dot{C}(t)V^2(t)dt,
\end{aligned}
\end{equation} 
which again leads to the passivity condition (\ref{eq:dotC_cond}). 

The dissipativity and passivity theory of the Appendix implies that for \emph{any} linear cyclicly varying capacitance $C(t)$, that is either described by (\ref{eq:LTVcap}), (\ref{eq:diff_VC}), or (\ref{eq:CpH}), there should \emph{always} exist a particular current excitation profile that let us extract more energy from the device than it is supplied with. This conclusion is completely in line with the one--port analysis of \cite{Georgiou2020}. An explicit construction of such a current function will be given next.

\begin{figure}[t]
\begin{center}
\includegraphics[width = 55mm]{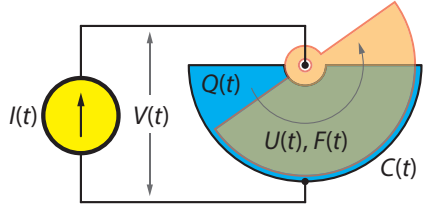}
\caption{A variable capacitor system.}
\label{fig:time_varying_cap_system}
\end{center}
\end{figure}

\subsection{An apparent `energy harvesting' device} 

Consider the variable capacitor system depicted in Figure~\ref{fig:time_varying_cap_system}. Two plates are shown, one of which is driven by a constant angular velocity $\omega = \frac{1}{2}$ [rad/s], and let us adopt the example in \cite[p.~304]{ChuaBook1987} in which the plate arrangement is assumed to be such that capacitance between the plates varies cyclicly according to 
\begin{equation}\label{eq:C(t)}
C(t) = 2 + \sin(\omega t).
\end{equation}

For a purely sinusoidal current $I(t) = \cos(2\omega t)$ and $Q(0) = 0$, the device models (\ref{eq:LTVcap}) or (\ref{eq:CpH}) appear to be lossless since the supplied energy (\ref{eq:LTVcap_suppliedE}) equals
\begin{equation*}
\oint\limits_{0}^{T} V(t)I(t) dt = 0,
\end{equation*}
where $T=2 \pi/\omega = 4 \pi$ [s] is the time of one full cycle of the variable capacitance. However, inspired by \cite{Jeltsema_Nature}, selecting the following (zero-mean) $4\pi$--periodic current profile
\begin{equation}\label{eq:I_profile}
\begin{aligned}
I(t) &= \frac{5}{2}\sin(\omega t) - \frac{1}{4}\cos(\omega t) -5\sin(2\omega t)\\ 
& \qquad - \frac{5}{4}\cos(2\omega t) + \frac{11}{4}\cos(3\omega t) - \frac{5}{4}\cos(4\omega t),
\end{aligned}
\end{equation}
yields that the supplied energy over one cycle equals
\begin{equation*}
\oint\limits_0^{T} V(t)I(t)dt \approx \boxed{-0.7 \ \text{[joules]}}
\end{equation*}
which shows that there is more energy gained from the device than it is supplied with. Moreover, for each subsequent cycle, the gained energy increases linearly with the number of cycles. This suggests that this device is an \emph{overunity} system, or, in a more classical parlance, a \emph{perpetual motion} machine of the first kind. 

\begin{figure}[h]
\begin{center}
\includegraphics[width = 80mm]{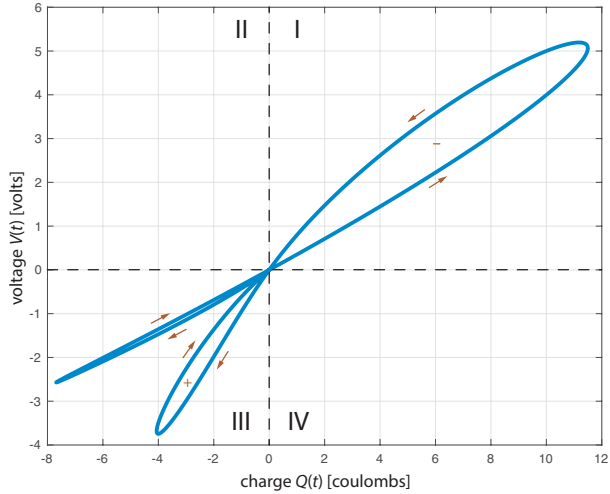}
\caption{Lissajous plot of the voltage versus the charge of the electromechanical system of Figure~\ref{fig:time_varying_cap_system} according to the (incomplete) device model (\ref{eq:LTVcap}) or (\ref{eq:CpH}).}
\label{fig:TVC_lissajous}
\end{center}
\end{figure}

The gain of energy is also evident from the Lissajous plot in Figure \ref{fig:TVC_lissajous}; the loop in the first quadrant has a counter--clockwise orientation (energy release), whereas the loops in the third quadrant both have a clockwise orientation (energy absorption); the energy release clearly exceeds the absorption of energy. Thus, we can apparently harvest `free' energy from a cyclicly changing capacitor device that is either described by the (incomplete) device model (\ref{eq:LTVcap}) or (\ref{eq:CpH}).   

\section{Now, get the physics right\ldots}\label{sec:completeness}

The obvious reason why an apparent amount of energy can be extracted from the cyclic capacitor device in Figure~\ref{fig:time_varying_cap_system} is that there is also some sort of mechanical torque required to change the dielectrics of the capacitor. However, this torque--like quantity is not taken into consideration in the device model (\ref{eq:LTVcap}) or (\ref{eq:CpH}). To fix this, the device model should be augmented with a mechanism that causes the change in the capacitance. This introduces a second degree of freedom that can be included by extending (\ref{eq:LTVcap}) or (\ref{eq:CpH}) with a state equation of the form
\begin{equation*}
\dot{C}(t) = U(t),
\end{equation*}
where $U(t)$ represents an independent external variable that causes the rate of change in the capacitance. 

\subsection{Time--varying capacitor as two--port device}

Concerning the device model (\ref{eq:CpH}), this will result in the following linear and time--invariant state equations\footnote{See \cite[p.~168]{passivitybook} for a similar type of system in the context of variable impedance control. }
\begin{subequations}\label{eq:TVC_complete_states}
\begin{empheq}[box=\fbox]{align}
\dot{Q}(t) & = I(t)\\
\dot{C}(t) & = U(t) \label{eq:dotC}
\end{empheq}
\end{subequations}

\medskip

\noindent with the nonlinear and time--invariant conjugated outputs (compare to (\ref{eq:CpH_dotS_full}))
\begin{equation}\label{eq:TVC_complete_outputs}
\boxed{V(t) = \frac{Q(t)}{C(t)}}  \ \ \text{and} \ \ \boxed{F(t) = -\frac{1}{2C^2(t)}Q^2(t)}
\end{equation}
and initial conditions $Q(0)=Q_0$ and $C(0)=C_0$, respectively. 

Note that $U(t)$ has the units of Farads per second [F/s], and as such can be interpreted as some kind of mechanical--like velocity. Consequently, the associated conjugate output $F(t)$ must represent a mechanical--like force or torque, known in classical electro--mechanics as the counter--electromotive force or back EMF. Indeed, in the light of (\ref{eq:TVC_complete_states})--(\ref{eq:TVC_complete_outputs}), the time--dependent stored energy (\ref{eq:TVC_S}) now matures to a proper nonlinear and time--invariant state function
\begin{equation}\label{eq:TVC_SQC}
S\big(Q(t),C(t)\big) = \frac{1}{2C(t)}Q^2(t),
\end{equation}
for which the associated power--balance becomes
\begin{equation}\label{eq:TVC_complete_S}
\frac{d}{dt}S\big(Q(t),C(t)\big) = V(t)I(t) + F(t)U(t),
\end{equation}
for all time $t$. 

The equality (\ref{eq:TVC_complete_S}) readily implies that the device model (\ref{eq:TVC_complete_states})--(\ref{eq:TVC_complete_outputs}) constitutes  a passive (in fact, \emph{lossless}) \emph{two--port} device, with conjugated electrical port variables $(I,V)$ and mechanical--like port variables $(U,F)$, and therefore the device model (\ref{eq:TVC_complete_states})--(\ref{eq:TVC_complete_outputs}) satisfies the First Law of thermodynamics---as should be expected. 

\subsection{Mechanical side of the picture}\label{sec:mech_side}

Starting from the constitutive relationship (\ref{eq:QvsV}), the two--port device model (\ref{eq:TVC_complete_states})--(\ref{eq:TVC_complete_outputs}) should be considered as the minimum set of equations that is needed to be energetically consistent. The physical structure of the mechanism that causes the capacitance to vary is still missing from the device model. Indeed, apart from assuming perfectly conducting electrical wires and frictionless mechanics, a practical rotating plate in Figure \ref{fig:time_varying_cap_system} will always experience some moment of inertia. As a result, there will also be storage of kinetic energy due to an induced angular momentum. 

Assuming the moment of inertia, say $J$, is constant, the angular momentum $P$ is given by the constitutive relationship $P(t) = J \dot{\Theta}(t)$,
where $\Theta(t)$ represents the angular displacement of the rotating plate. Consequently, the capacitance can be expressed as an explicit function of the angular displacement, so that (\ref{eq:QvsV}) is replaced by a constitutive relationship of the form
\begin{equation}
Q(t) = C\big(\Theta(t)\big)V(t)
\end{equation}
and (\ref{eq:dotC}) is replaced by the dynamics of the mechanical subsystem. 

For details and other examples that involve state--modulated capacitors, such as capacitor microphones and micro electro--mechanical systems (MEMS), the interested reader is referred to \cite{passivitybook}  and \cite{TAC2022}, and the references therein.

\section{A two capacitor--like `paradox'}\label{sec:2Cap_paradox}

\begin{figure}[h]
\begin{center}
\includegraphics[width = 50mm]{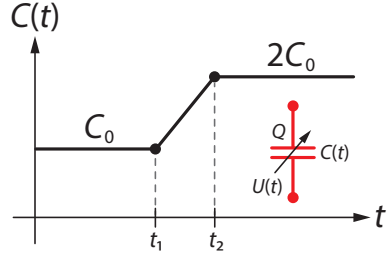}
\caption{A capacitor with a linear change in its dielectrics.}
\label{fig:C_2_2C}
\end{center}
\end{figure}

Although we mainly concentrated on capacitors that vary in a continuous cyclic fashion, the proposed two--port device model (\ref{eq:TVC_complete_states})--(\ref{eq:TVC_complete_outputs}) is also valid for discontinuous and non--cyclic capacitive devices. To illustrate this, consider a capacitor that is initially charged with a certain amount of charge, say $Q$, and at a certain moment in time, say $t_1 = 0$, changes linearly to twice its capacitance at $t_2 = 1$; see Figure \ref{fig:C_2_2C}. 

Since the total charge is conserved, the dynamics covering this time--varying process is
\begin{subequations}\label{eq:changing_C}
\begin{align}
\dot{Q}(t) &= 0,\\
\dot{C}(t) &= U(t), 
\end{align}
\end{subequations} 
with external `mechanical' input
\begin{equation}
U(t) = \left\{
\begin{array}{ccc}
0 & \text{for} & t < 0,\\ 
1/T & \text{for} & 0 \leq t < T,\\
0 & \text{for} & t \geq T,
\end{array}\right.
\end{equation}
with $T = t_2 - t_1$.  
If the capacitance at $t \leq t_1$ is $C_0$, then the energy stored in the capacitor at $t = t_1$ equals 
\begin{equation}\label{eq:TVC_S_t1}
S\big(C(t_1)\big) = \frac{1}{2C_0}Q^2.
\end{equation} 
However, after the transition from $C$ to $2C$ is completed, the energy is reduced to
\begin{equation}\label{eq:TVC_S_t2}
S\big(C(t_2)\big) = \frac{1}{4C_0}Q^2.
\end{equation}
Now, where did this energy go? 

Again, the most natural answer lies in the fact that we did not take the mechanism that causes the change of the capacitance into account. Indeed, differentiating the total energy at any arbitrary time $t$ along the trajectories of (\ref{eq:changing_C}) yields
\begin{equation}
\frac{d}{dt}\left(\frac{1}{2C(t)}Q^2\right) =  -\frac{1}{2}\frac{Q^2}{C^2(t)}\dot{C}(t) \equiv F(t)U(t),
\end{equation}
which, in turn, implies that
\begin{equation}
S\big(C(t_2)\big) - S\big(C(t_1)\big) = \int\limits_{t_1}^{t_2} F(t)U(t)dt = -\frac{1}{4C_0}Q^2. 
\end{equation}
Interestingly, note that in the limit $t_2 \to t_1$, and subsequently, 
\begin{equation}
\lim_{T \to 0} U(t) = \delta(t),
\end{equation}
with $\delta(t)$ the Dirac impulse function, the result still holds and reduces to a similar problem known as the ``Two Capacitor Paradox;'' see \cite[p.~155]{Hambley} and \cite{Hoekstra} and the references therein. Hence, the resolution to this paradox is that the apparent `lost' energy is due to the mechanism that causes the change in the system. This is to be expected, and in full agreement with the analysis of \cite{Hoekstra}, since `mechanical' work is required to change the dielectrics of the capacitor or to operate an ideal switching device that connects to identical capacitors $C_0$, one of which is initially charged to a certain charge level $Q$ and the other initially uncharged, at a certain instant of time.

\section{Final remarks}

The main conclusion to be drawn from the perspectives of this paper is that the classical charge--voltage relationship (\ref{eq:QvsV}) is generally valid---at least within the context of circuit--theoretic device models and their postulates. However, a naive application of the one--port device model (\ref{eq:LTVcap}) generally leads to incomplete and thermodynamically inconsistent results. To avoid such issues, the proposed two--port model (\ref{eq:TVC_complete_states})--(\ref{eq:TVC_complete_outputs}) should be used instead. 

On the other hand, as briefly outlined in Section \ref{sec:mech_side}, if one includes physical details about the mechanism that causes the change in capacitance, then a capacitor that explicitly depends on time generally loses its relevance.  

Furthermore, although only time--varying capacitance is discussed, the same analysis can be carried over verbatim to a time--varying inductance that is described by a constitutive relationship of the form $\Phi(t) = L(t)I(t)$ and a (naive) application of Faraday's law, $V(t) = \dot{\Phi}(t)$, yielding
\begin{equation}
V(t) = \frac{d}{dt}\big[L(t)I(t)\big] = L(t)\dot{I}(t) + \dot{L}(t)I(t),
\end{equation}
where $\Phi(t)$ denotes the flux--linkage. 

Finally, concerning the constitutive relationship (\ref{eq:Q*V}) as proposed in \cite{Vikash2022}, it is well--known, see e.g.,~\cite{DeHoop}, that in electromagnetic media there are situations that the material itself may have some sort of inherent memory. In such case, the constitutive relationships associated to the permittivity, permeability and conductivity might involve convolutions between the intensive and extensive field variables. However, it is rather doubtful if such cases can be simplified to quasi--static circuit--theoretic models without violating any physical laws.

\section*{Appendix: Dissipativity and passivity theory}\label{appendix}

Consider a system with state vector $x$ and a vector of external (e.g., input and output) variables $w$. Consider a scalar--valued {\it supply rate} $s(w)$. A function $S(x)$ is said to be a {\it storage function} (with respect to the supply rate $s$) if along all trajectories of the system and for all $t_1 \leq t_2$ and $x(t_1)$ it satisfies the {\it dissipation inequality} (see e.g., \cite{hillmoylan1980}, \cite{passivitybook}, and \cite{willems1972})
\begin{equation}\label{diss}
S\big(x(t_2)\big) - S\big(x(t_1)\big) \leq \int\limits_{t_1}^{t_2} s\big(w(t)\big) dt.
\end{equation}
Interpreting $s\big(w(t)\big)$ as `power' supplied to the system at time $t$, and $S\big(x(t)\big)$ as stored `energy' while at state $x(t)$, this means that increase of the stored energy can only occur due to externally supplied power. 

Following \cite{willems1972}. the system is called {\it dissipative} (with respect to the supply rate $s$) if there exists a {\it nonnegative} storage function $S$.\footnote{Since addition of an arbitrary constant to a storage function again leads to a storage function, the requirement of nonnegativity of $S$ can be relaxed to $S$ being {\it bounded from below}.} Furthermore, it is called {\it lossless} (with respect to $s$) if there exists a nonnegative storage function satisfying the dissipation inequality \eqref{diss} with {\it equality}.

An external characterization of dissipativity, in terms of the behavior of the $w$ trajectories, is the following \cite{willems1972}. Define for any $x$ the expression (possibly infinite)
\begin{equation}
S_a(x):= \sup_{w, T\geq 0} - \int\limits_0^T s\big(w(t)\big) dt,
\end{equation}
where the supremum is taken over all external trajectories $w(\cdot)$ of the system corresponding to initial condition $x(0)=x$, and all $T\geq 0$. Obviously, $S_a(x)\geq 0$. Then the system is dissipative with respect to the supply rate $s$ if and only if $S_a(x)$ is {\it finite} for every $x$. 

Interpreting as above $s(w)$ as `power' supplied to the system, $S_a(x)$ is the maximal `energy' that can be extracted from the system at initial condition $x$, and  the system is dissipative if and only if this maximal `energy' is finite for any $x$. Furthermore, if  $S_a(x)$ is finite for every $x$ then $S_a$ is itself a non-negative storage function, and is in fact the {\it minimal} nonnegative storage function (generally among many others). 

Dropping the requirement of non-negativity of the storage function $S$ leads to the notion of {\it cyclo--dissipativity}, respectively {\it cyclo--losslessness} \cite{passivitybook, willems1973}. First note that if a general, possibly indefinite, storage function $S$ satisfies \eqref{diss}, then for all {\it cyclic} trajectories, i.e., such that $x(t_1)=x(t_2)$, we have that
\begin{equation}\label{cyclic}
\oint\limits_{t_1}^{t_2}  s\big(w(t)\big) dt \geq 0,
\end{equation}
which holds with equality in case \eqref{diss} holds with equality. This leads to the following external characterization of {\it cyclo}-dissipativity and {\it cyclo}--losslessness.

A system is {\it cyclo-dissipative} if (\ref{cyclic}) holds
for all $t_2\geq t_1$ and all external trajectories $w$ such that $x(t_2)=x(t_1)$. In case \eqref{cyclic} holds with equality, we speak about cyclo--losslessness. Furthermore, the system is called cyclo-dissipative {\it with respect to} $x^*$ if \eqref{cyclic} holds for all $t_2\geq t_1$ and all external trajectories $w$ such that $x(t_2)=x(t_1)=x^*$, and cyclo--lossless with respect to $x^*$ if this holds with equality. 

Interpreting again $s(w)$ as the power provided to the system, cyclo--dissipativity thus means that for any cyclic trajectory the net amount of energy supplied to the system is non-negative, and zero in case of cyclo--losslessness.

In case of the {\it passivity} supply rate 
\begin{equation}
s(w)=s(u,y)=y^Tu, 
\end{equation}
where $w=(u,y)$ and $u$ and $y$ are equally dimensioned vectors, the terminology `dissipativity' in all of the above is replaced by the classical terminology of {\it passivity}. In this case $u$ and $y$ typically are vectors of power--conjugate variables, like forces and velocities, and voltages and currents. 


\begin{thebibliography}{00}
\bibitem{Allagui2022} A. Allagui, A.S. Elwakil, C. Psychalinos,
Comment on ``Origin of the Curie?von Schweidler law and the fractional capacitor from time-varying capacitance'' [J. Pow. Sources 532 (2022) 231309], Journal of Power Sources, Volume 551, 2022, 232166.

\bibitem{Biolek} D.  Biolek, Z. Kolka, V. Biolkova, ``Modeling time-varying storage components in PSpice,'' In Proc. Electronic Devices and Systems IMAPS CS International Conference EDS, pp.~39--44, 2007. 

\bibitem{ChuaBook1987} L.O. Chua, C.A. Desoer and E.S. Kuh, {\it Linear and Nonlinear Circuits}, McGraw--Hill Inc., 1987. 

\bibitem{Desoer} C.A. Desoer and E.S. Kuh, {\it Basic Circuit Theory}, McGraw--Hill Book Company, 1969.

\bibitem{Fouda2020} M.E. Fouda, A. Allagui, A.S. Elwakil, S. Das, C. Psychalinos, and A.G. Radwan, ``Nonlinear charge--voltage relationship in constant phase element,'' {\it AEU-Int. J. Electron. Commun.}, 117:153104, 2020.

\bibitem{Georgiou2020} T.T. Georgiou, F. Jabbari and M.C. Smith, ``Principles of lossless adjustable one--ports,'' {\it IEEE Transactions on Automatic Control}, vol. 65, no. 1, pp. 252-262, Jan. 2020. 

\bibitem{Hambley} A.R. Hambley, {\it Electrical Engineering: Principles and Applications}, 6th Edition, Pearson Education Limited, 2014. 

\bibitem{hillmoylan1980}
D.J. Hill and P.J. Moylan,
``Dissipative dynamical systems: basic input-output and state properties,''
{\it Journal of the Franklin Institute}, 309(5):327-357, 1980.

\bibitem{Hoekstra} J. Hoekstra, ``A solution of the two--capacitor problem through its similarity to single-electron electronics,'' {\it IEEE Open Journal of Circuits and Systems}, Vol. 1, pp. 13--21, 2020, DOI: 10.1109/OJCAS.2020.2977216.

\bibitem{DeHoop} A.T. de Hoop, {\it Handbook of Radiation and Scattering of Waves}, Academic Press, 1995. 

\bibitem{Jeltsema2009} D. Jeltsema and J.M.A. Scherpen, ``Multidomain modeling of nonlinear networks and systems,'' {\it IEEE Control Systems Magazine}, Vol. 29, No. 4, pp. 28--59, July 2009.  

\bibitem{Jeltsema_Nature} D. Jeltsema and A.J. van der Schaft, ``Ideal memcapacitors and meminductors are overunity devices,'' {\it Scientific Reports}, 10, 16688, 2020, DOI: 10.1038/s41598--020--73833--3.

\bibitem{Vikash2022} V. Pandley, ``Origin of the Curie--von Schweidler law and the fractional capacitor from time--varying capacitance, {\it Journal of Power Sources},  532, 231309, 2022. 

\bibitem{Vikash2022_response} V. Pandey, Response to ``Comment on `Origin of the Curie--von Schweidler law and the fractional capacitor from time--varying capacitance [J. Power Sources 532 (2022) 231309]' '', Journal of Power Sources, Volume 551, 2022, 232167. 

\bibitem{Perkins1982} D.H. Perkins, {\it Introduction to High-Energy Physics}, Addison-Wesley, 1982. 

\bibitem{passivitybook}
A.J. van der Schaft, {\it $L_2$-Gain and Passivity Techniques in
Nonlinear Control}, 3rd Edition, Springer International, 2017.

\bibitem{NOW_book}
A.J. van der Schaft and D. Jeltsema, {\it Port-Hamiltonian Systems Theory: An Introductory Overview}, Foundations and Trends in Systems and Control,
NOW Publishers, 2014. 

\bibitem{TAC2022} A.J. van der Schaft and D. Jeltsema, ``Limits to Energy Conversion,'' in {\it IEEE Trans. on Automatic Control}, vol. 67, no. 1, pp. 532--538, Jan. 2022, DOI: 10.1109/TAC.2021.3075652.

\bibitem{Stewart} J. Stewart, {\it Calculus}, 6th Edition, Thomson Brooks/Cole, 2008. 

\bibitem{Valkenburg1955} M.E. van Valkenburg, {\it Network Analysis}, Prentice--Hall, Inc., 1955. 

\bibitem{willems1972}
J.C. Willems, ``Dissipative dynamical systems, part I: general theory,'' {\it Arch. Rat. Mech. and Analysis}, 45(5):321--351, 1972.

\bibitem{willems1973}
J.C. Willems, ``Qualitative behavior of interconnected systems,''
{\it Annals of Systems Research}, 3:61--80, 1973.
\end{thebibliography}
\end{document}